\documentclass[10pt]{article} 
\usepackage{amscd} 
\usepackage{amssymb} 
\usepackage{amsmath} 
\usepackage{amsthm} 
\usepackage{amsbsy} 
\usepackage{graphics} 

\newcommand{\R}{\ensuremath{\mathbb R}}
 
\newcommand{\Z}{\ensuremath{\mathbb Z}} 
 
\newtheorem{theorem}{Theorem}[section] 
\newtheorem{proposition}{Proposition}[section]  
\newtheorem{corollary}{Corollary}[section]  
  
\newtheorem{remark}{Remark}[section]  
\newtheorem{lemma}{Lemma}[section]  
 
\newtheorem{problem}{Problem} 
\newtheorem{conjecture}{Conjecture}

\begin{document} 
\title{Surface cubications mod flips}  
\author{Louis Funar
\\ 
\small \em Institut Fourier BP 74, UMR 5582 CNRS,\\ 
\small \em University of Grenoble I,\\ \small \em 
38402 Saint-Martin-d'H\`eres cedex, France \\ \small \em e-mail: {\tt 
funar@fourier.ujf-grenoble.fr} \\ 
} 
{\small\date{\today}} 
 
\maketitle 
{\abstract Let $\Sigma$ be a compact surface.  
We prove that the set of  marked surface cubications modulo flips, up to  
isotopy, is in one-to-one correspondence  with  
$\Z/2\Z\oplus H_1(\Sigma,\partial \Sigma;\Z/2\Z)$. }

\vspace{0.3cm} 
 
\vspace{0.3cm} 

\section{Introduction and statements}
{\bf Cubical complexes and marked cubications.} 
A {\it cubical complex} is a  finite dimensional complex  
$C$ consisting of Euclidean cubes, 
such that the intersection of  two of its  cubes  is a finite union of cubes 
from  $C$, once a cube is in $C$ then all its faces belong to $C$ and 
each point  has a neighborhood intersecting only finitely many cubes of $C$.   
A {\it cubication} of a  topological manifold is a cubical 
complex  that is homeomorphic to the manifold. If the manifold is a PL
manifold then one requires that the cubication  be combinatorial  and 
compatible with the PL structure. Our definition of cubication 
is slightly more general than the usual one, because we don't  
require that the intersection of two cubes consists of a single cube but only  
a {\em finite union} of cubes.

\vspace{0.2cm}\noindent
The study of simplicial complexes and manifold triangulations 
lay at the core of combinatorial topology. Cubical 
complexes and cubications might offer an alternative approach since, despite 
their similarities,  they present some new features.  

\vspace{0.2cm}\noindent
Any triangulated manifold admits a cubication, since we can decompose 
an $n$-dimensional simplex $\Delta^n$ into $n+1$ cubes of dimension $n$. 
For $k=1$ to $n$ we adjoin, inductively, the barycenter of each $k$-simplex in 
$\Delta^n$ and join it with the barycenters of its 
faces. This way we obtain the 1-skeleton of a cubical complex, as shown in the 
figure below for $n=3$. 

\vspace{0.2cm}
\begin{center}
\includegraphics{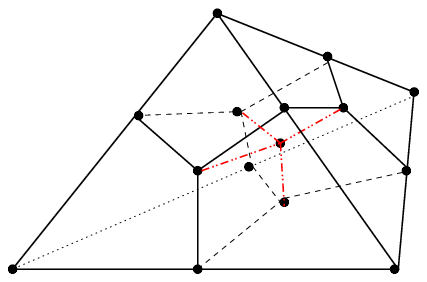}
\end{center}

\vspace{0.2cm}\noindent
Thus,  roughly speaking, 
working with simplicial complexes is equivalent to working with 
cubical complexes, from topological viewpoint. However, 
we will show that cubications 
encode additional topological information.

\vspace{0.2cm}\noindent
It will be more convenient in the sequel to  work with 
marked cubications instead of cubications. 
A {\em marked} cubication of the manifold $M$ 
consists of a couple $(C,\varphi)$, where
$C$ is a cubication and $\varphi:|C|\longrightarrow M$ is 
a PL homeomorphism  (called the {\em marking}) of 
its subjacent space $|C|$ onto $M$. 
The marked cubications $(C,\varphi)$ and $(C',\varphi')$ are said to be 
{\em isotopic} if there exists a combinatorial isomorphism 
$j:C\to C'$ between the two cubical complexes and a PL homeomorphism 
$\Phi$ of $M$ 
such that $\Phi\circ \varphi=\varphi'\circ J$, 
where $J:|C|\to |C'|$ is the PL homeomorphism induced by $j$, and 
both $J$ and $\Phi$ are isotopic to identity.
The isotopy class of the image by $\varphi$ of the  
skeleton of $C$ in $M$ determines the isotopy class 
of the marked cubulation $(C,\varphi)$.  
Thus, marked cubications underlying a given cubication $C$ are acted upon  
transitively by the mapping class group of $M$.

\vspace{0.2cm}\noindent
{\bf Bi-stellar moves.} 
We will consider below PL manifolds, i.e.  topological manifolds  
endowed with triangulations (called combinatorial) 
for which the link of each  vertex is PL homeomorphic to the 
boundary of the simplex. Recall that two simplicial complexes 
are PL homeomorphic if they admit combinatorially isomorphic subdivisions. 
There exist topological manifolds which have several 
PL structures and it is still unknown 
whether all topological manifolds have triangulations (i.e. whether they 
are homeomorphic to simplicial complexes), without requiring 
them to be combinatorial. 

\vspace{0.2cm}\noindent
It is not easy to decide whether two given triangulations 
define or not the same  PL structure.  One difficulty is that one has to 
work with arbitrary subdivisions and there are 
infinitely many distinct combinatorial types of such.  
In the early sixties one looked upon a more convenient 
set of transformations permitting to connect  PL equivalent 
triangulations of a given manifold. The  simplest proposal 
was the so-called  {\em bi-stellar} moves  which are defined for 
$n$-dimensional complexes, as follows: we excise $B$ and replace it 
by $B'$, where $B$ and $B'$ are complementary balls that are unions of 
simplexes in the boundary $\partial \Delta^{n+1}$ of the standard  
 $(n+1)$-simplex.  It is obvious that such transformations 
do not change the  PL homeomorphism type of the complex. 
Moreover, U.Pachner (\cite{pa1,pa2}) proved in 1990  that conversely, 
any two  PL triangulations of a PL manifold (i.e. 
the two triangulations define the same PL structure) can be connected by a 
sequence of bi-stellar moves. 
One far reaching application of Pachner's theorem was the 
construction of the Turaev-Viro quantum invariants  (see \cite{Tu}) 
for 3-manifolds.

\vspace{0.2cm}\noindent 
{\bf Habegger's problem on cubical decompositions.}
It is natural to wonder whether a similar result holds  
for cubical decompositions, as well. The cubical decompositions 
that we consider will be PL decompositions that define the same 
PL structure of the manifold. 

\vspace{0.2cm}\noindent 
Specifically, N.Habegger asked (\cite{Kir}, problem 5.13) the following: 
\begin{problem}
{\em Suppose  that we have two PL cubications of the same PL manifold. Are 
they related by the following set of moves: excise $B$ and replace it 
by  $B'$, where $B$ and $B'$ are complementary balls (union of $n$-cubes) 
in the boundary of the standard $(n+1)$-cube?}
\end{problem} 
 
\vspace{0.2cm}\noindent 
These moves have been called {\em cubical} or {\em bubble}
moves in \cite{F,F2}, and  (cubical) {\em flips} in \cite{BEE}.
Notice that the flips did already appear in the 
mathematical polytope literature (\cite{BS,Z}). 

\vspace{0.2cm}\noindent  
The problem above was addressed in (\cite{F,F2}), where we show that,  
in general, there are  topological obstructions for 
two cubications being flip equivalent.

\vspace{0.2cm}\noindent  
Notice that  acting by cubical flips one can create  cubications 
where cubes have several faces in common or pairs of faces of the same 
cube are identified. Thus we are forced to allow this greater degree 
of generality in our definition of cubical complexes.

\vspace{0.2cm}\noindent
{\bf Related work on cubications.} 
In the meantime this and related problems have been approached by 
several people working in computer science or combinatorics 
of polytopes (see \cite{BEE,CMP,E,Mi,SZ}). Notice also that 
the 2-dimensional case of the sphere $S^2$ was actually solved earlier 
by Thurston (see \cite{Th}).  
Observe that there are several terms in 
the literature describing the same object. For instance 
the cubical decompositions of surfaces are also called {\em quadrangulations} 
(\cite{Na1,Na2,NO}) or  
{\em quad  surface meshes}, while 3-dimensional cubical complexes are called 
{\em hex meshes} in the computer science papers (e.g. \cite{BEE}). 
We used the term {\em cubulation} in \cite{F,F2}.

\vspace{0.2cm}\noindent 
Remark that there is some related work that has been done by Nakamoto 
 (see  \cite{Na1,Na2,NO}) 
concerning the equivalence of cubications 
{\em of the same order } by means of two transformations (that preserve 
the number of vertices):   the {\em diagonal slide},  
in which one exchanges one diameter of a hexagon for another, and 
the  {\em diagonal rotation},  
in which the neighbors of a vertex of degree 2 inside 
a quadrilateral are switched. In particular, it was proved that 
any two  cubications of a closed  orientable surface can be transformed 
into each other, up to isotopy, by diagonal slides and diagonal rotations 
if they have the same (and sufficiently large) number of vertices 
and if their 1-skeleta define the same mod 2 homology classes. 
Moreover, one can do 
this while preserving the simplicity of the cubication (i.e.
not allowing double edges).

\vspace{0.2cm}\noindent 
{\bf Immersions and cobordisms.}
Let $M$ be a $n$-dimensional  manifold. Consider the set of 
immersions $f:F\rightarrow M$ with $F$ a closed $(n-1)$-manifold. 
Impose on it the following equivalence relation: $(F,f)$ is {\em cobordant} to 
$(F',f')$ if there exist 
 a cobordism $X$ between $F$ and $F'$, that is, a 
compact $n$-manifold $X$ with boundary $F\sqcup F'$, and an immersion $\Phi:X\rightarrow M\times I$, transverse to 
the boundary, 
 such that $\Phi_{|F}=f\times\{0\}$ 
and $\Phi_{|F'}=f'\times\{1\}$. 
 
\vspace{0.2cm}\noindent 
Once the manifold $M$ is fixed, the 
set $N(M)$ of cobordism classes of codimension-one immersions  in $M$ is 
an abelian group with the composition law given by disjoint union.

\vspace{0.2cm}\noindent 
{\bf Cubications vs immersions' conjecture.}
Our approach in \cite{F} to the flip equivalence problem  
aimed at finding a general solution in terms of  some algebraic 
topological invariants.   Specifically,  
we stated (and proved half of) the following  conjecture: 

\begin{conjecture}
The set of  marked cubical decompositions of the  closed manifold $M^n$
modulo cubical flips is in bijection with the elements of the cobordism group 
of codimension one immersions into $M^n$. 
\end{conjecture}

\vspace{0.2cm}\noindent 
The solution of this conjecture would lead to a 
quite satisfactory answer to the problem of Habegger.  

\vspace{0.2cm}\noindent 
Notice that, when a cubical move is performed on the cubication $C$ endowed 
with a marking, there is a natural marking induced for the flipped cubication. 
Thus it makes sense to consider the set of marked cubications mod flips.

\vspace{0.2cm}\noindent 
We proved in \cite{F} the existence of a surjective map between 
the two sets.

\vspace{0.2cm}\noindent 
{\bf Digression on smooth vs PL category.}  
There might be several possible interpretations for   
the conjecture above. We can work, for instance, 
in  the PL category and thus cubications, 
immersions and cobordisms are supposed PL. Generally speaking   
there is little known about the cobordism group of 
PL immersions and their associated Thom spaces, in comparison 
with the large literature on smooth immersions. However, in 
the specific case of codimension-one immersions we are able to 
compare the relevant bordism groups. If $M$ is a compact $n$-dimensional 
manifold then the bordism group $N_k(M)^{\rm PL}$ 
of $PL$ codimension-$k$ immersions up to PL cobordisms 
is given  in homotopy theoretical terms by the formula: 
\[   N_k(M)^{\rm PL}= [M, \Omega^{\infty}S^{\infty} MPL(k)] \]
where $PL(k)$ is the semi-simplicial group of PL germs of maps on $\R^k$, 
$MPL(k)$ is a suitable Thom space associated to it (see e.g. 
\cite{Wil}), $\Omega$ denotes the loop space and $S$ denotes the 
reduced suspension,
while $[X,Y]$ denotes the set of homotopy classes 
of maps $X\to Y$. 
This follows along the same lines as the results of Wells 
(\cite{WelCGI}), where is considered only the smooth case, by using instead 
of the classical Smale-Hirsch theory on smooth immersions the 
Haefliger-Poenaru classification of combinatorial immersions  from 
(\cite{HP}). 

\vspace{0.2cm}\noindent 
On the other hand, when $M$ is smooth, the bordism group $N_k(M)$ of 
codimension-$k$ smooth immersions is given by the similar formula 
from (\cite{WelCGI}): 
\[   N_k(M)= [M, \Omega^{\infty}S^{\infty} MO(k)] \]
where $MO(k)$ is the Thom space associated to the orthogonal group $O(k)$. 

\vspace{0.2cm}\noindent 
From the general results of Kuiper and Lashof (\cite{KuLa}) concerning the 
unstable homotopy type of $PL(k)$ one obtains that the natural 
inclusion map $O(1)\hookrightarrow PL(1)$ induces a weak homotopy 
equivalence $MO(1)\to MPL(1)$. This result was improved later 
by Akiba, Scott and Morlet (\cite{A,Mo,Sc}) to weak homotopy 
equivalences $MO(k)\to MPL(k)$ for all $k\leq 3$. 

\vspace{0.2cm}\noindent 
This shows that codimension-one immersions of PL manifolds into 
a smooth manifold are PL cobordant to a smooth immersion and also 
that the existence 
of a PL cobordism between two  smooth immersions implies the existence 
of a smooth cobordism.  
Consequently, when the manifold $M$ is smooth, we can use either PL or smooth 
immersions and bordisms, as the associated groups are naturally isomorphic.

\vspace{0.2cm}\noindent 
{\bf Smooth cubications.}
However, in the DIFF category it is appropriate to consider only  
those cubications which are smooth. 
Smooth cubications are defined following Whitehead's definition of 
smooth triangulations  from (\cite{Wh1}), but we have to change it 
slightly in order to apply to the more general cubical complexes considered 
here. 

\vspace{0.2cm}\noindent 
Let $M$ be a smooth manifold and $C$ a cubical complex. 
A map $f:|C|\to M$  is  called {\em smooth} if the restriction  of $f$ to each 
cube of $C$ is smooth. Moreover, $f$ is {\em non-degenerate} 
if all these restrictions are of maximal 
rank. Finally $f:|C|\to M$ is a {\em smooth cubication} of $M$ if $f$ is a 
non-degenerate  homeomorphism onto $M$. According to (\cite{Mu}, Theorem 8.4) 
this definition is equivalent to Whitehead's one, when applied to simplicial 
complexes.

\vspace{0.2cm}\noindent 
In small dimensions (e.g. when the dimension is at most 3) the PL and 
DIFF categories are equivalent. In particular, we can 
assume from now on that we are working in the DIFF category and all 
objects are smooth, unless the opposite is explicitly stated.

\vspace{0.2cm}\noindent 
{\bf Computations of the cobordism group of immersions.}
Finding the cobordism group $N(M^n)$ 
of  (smooth) codimension-one immersions into the $n$-manifold $M^n$ 
was reduced to a homotopy problem by the results of 
\cite{VogCI,WelCGI}, as explained above. 
However, these techniques seem awkward to apply when one is looking for 
effective results. 
The group  $N(S^n)$ of codimension-one 
immersions in the $n$-sphere, up to cobordism,  is the 
 the $n$-th stable  homotopy group of ${\mathbb R\mathbb P}^{\infty}$ 
(since the Thom space $MO(1)$ is homotopy equivalent to 
${\mathbb R\mathbb P}^{\infty}$) 
and it was computed by Liulevicius (\cite{LiuTHA}) 
for $n\leq 9$ as follows:

\vspace{0.5cm} 
\begin{center} 
\begin{tabular}{|c|c|c|c|c|c|c|c|c|c|}\hline 
$n$ &  1 & 2 & 3 & 4 & 5 & 6 & 7 & 8 & 9\\ \hline 
$N(S^n)$ &$ {\Z}/2{\Z}$ & ${\Z}/2{\Z}$ & ${\Z}/8{\Z}$ 
&${\Z}/2{\Z} $ &$ 0$ & ${\Z}/2{\Z}$ & ${\Z}/16{\Z}\oplus {\Z}/2{\Z} 
$ &$ ({\Z}/2{\Z})^{\oplus 3}$ & $({\Z}/2{\Z})^{\oplus 4} $\\ \hline 
\end{tabular} 
\end{center} 
\vspace{0.5cm}

\vspace{0.2cm}\noindent 
It is known that, if $M^2$ denotes a closed surface, then: 
\[ N(M^2)\cong H_1(M^2, {\Z}/2{\Z})\oplus H_2(M^2, {\Z}/2{\Z})\]
Using geometric methods Benedetti and Silhol 
(\cite{BeSSPS}) and further Gini (\cite{Gin}) proved that, if $M^3$ is a 3-manifold, then 
\[N(M^3)\cong H_1(M^3, {\Z}/2{\Z}) \oplus H_2(M^3, {\Z}/2{\Z}) 
\oplus 
H_3(M^3,{\Z}/8{\Z})\]
the right side groups being endowed with a twisted product. 
The result has been extended to higher dimensional manifolds 
in \cite{FG}. 

\vspace{0.2cm}\noindent 
{\bf Manifolds with boundary.}
Habegger's problem from above makes sense also for 
PL cubications of manifolds with boundary. The question is whether 
two cubications that induce the same cubication on the 
boundary are flip equivalent. 

\vspace{0.2cm}\noindent 
Let $M$ be a compact $n$-dimensional  manifold with boundary $\partial M$. 
We will consider then the {\em proper}
immersions $f:F\rightarrow M$ with $F$ a compact $(n-1)$-manifold with 
boundary $\partial F$. This means that $\partial M$ is transversal 
to $f$ and $f^{-1}(\partial M \cap f(F))=\partial F$.  

\vspace{0.2cm}\noindent In order to define the cobordism equivalence for 
proper immersions we need to introduce more 
general immersions and manifolds. The 
compact $n$-manifold $X$ is a manifold with {\em corners} 
if $X$ is a PL manifold whose boundary $\partial X$ has a splitting 
$\partial X=F \cup \partial F \times [0,1] \cup F'$, where 
$\partial F=\partial F'$, and $F,F'$ are manifolds with boundary. 
One says that  
$F\cup F'$ is the {\em horizontal boundary} $\partial_HX$, 
$\partial F \times [0,1]$ 
is the {\em vertical boundary} $\partial_VX$  and their intersection
$\partial F \times \{0,1\}$, is the {\em corners} set.  

\vspace{0.2cm}\noindent 
Observe now that $M\times [0,1]$ is naturally a manifold with corners 
if $M$ has boundary, by using the splitting 
$\partial (M\times [0,1])=M \times \{0\} \cup \partial M\times [0,1] 
\cup M\times \{1\}$. Let $X$ be as above. One defines then an immersion 
$\Phi:X\to M\times [0,1]$ to be an immersion of manifolds with 
corners if 
\begin{enumerate}
\item $\Phi$ is proper and preserves the boundary type, by sending 
the horizontal (resp. vertical)  part into the horizontal (resp. vertical) 
boundary. Moreover, $\Phi$ is transversal to the boundary.  
\item  The restriction to the vertical part 
$\Phi:\partial F\times [0,1] \to \partial M\times [0,1]$ is a product
i.e. it is of the form 
$\Phi(x,t)=(\Phi(x,0),t)$. 
\end{enumerate}

\vspace{0.2cm}\noindent 
We say that  the proper immersion $(F,f)$ is {\em cobordant} to 
$(F',f')$ if  there exist a cobordism $X$ between $F$ and $F'$ 
(and thus $\partial F=\partial F'$)  which is a manifold with corners 
and a proper immersion of manifolds with corners 
$\Phi:X\rightarrow M\times I$, such that $\Phi_{|F}=f\times\{0\}$ 
and $\Phi_{|F'}=f'\times\{1\}$.

\vspace{0.2cm}\noindent 
The set of cobordism classes of immersions of codimension 
one  manifolds with boundary into a given manifold $M^n$  
with prescribed boundary immersion can be computed using the methods 
of \cite{FG}.

\vspace{0.2cm}\noindent 
{\bf The main result.}
The aim of this paper is to solve the extension 
of the cubications vs immersions conjecture 
in the case of compact surfaces, possibly with boundary.

\begin{theorem}
The set of  marked cubications of  the compact surface $\Sigma$ 
with prescribed boundary 
mod cubical flips is in one to one correspondence with the elements of 
$\Z/2\Z\oplus H_1(\Sigma, \partial \Sigma;\Z/2\Z)$. 
\end{theorem}

\vspace{0.2cm}\noindent 
The proof of this theorem, although elementary, uses  some methods 
from geometric topology and Morse theory.

\begin{remark}
One can identify a marked cubication with an embedding of a 
connected graph  in the surface, whose complementary is made of 
squares. The theorem says that any two graphs like that are related by 
a sequence of cubical flips and an isotopy of the surface.  
\end{remark}

\vspace{0.2cm}\noindent 
{\bf Acknowledgements.} We are grateful to the referee for many valuable 
comments and suggestions leading to a better presentation and 
the simplification of some proofs.


\section{Outline of the Proof}
{\bf Immersions associated to cubications}. We associate  
to each marked cubication $C$ of the $n$-dimensional manifold 
$M$ a codimension-one  
generic immersion $\varphi_C: N_C\longrightarrow M$ (the cubical complex $N_C$ 
is also called the derivative complex in \cite{BC}) of a manifold $N_C$ 
having one dimension less  than $M$. Here is the construction.   
Each cube is divided into $2^n$ equal cubes by $n$ hyperplanes 
which we call sections. When gluing together cubes in a cubical 
complex the sections are glued accordingly.  Then the  union of 
the  hyperplane sections  form the image $\varphi_C(N_C)$ of a 
codimension-one generic immersion. 
The cubulated manifold
$N_C$ is constructed as follows: consider the disjoint union of 
a set of $(n-1)$-cubes which is in bijection with the set of 
all sections, then glue together two $(n-1)$-cubes if their 
corresponding sections are adjacent in $M$. The immersion $\varphi_C$
is tautological: it sends a cube of  $N_C$ into the corresponding
section.  If the cubication $C$ is smooth then $N_C$ 
has a smooth structure and the immersion $\varphi_C$ can be  made smooth
by means of a small isotopy. 
This connection between cubications and immersions 
appeared independently in \cite{BC,F} but this was presumably known 
to specialists long time ago (see e.g. \cite{Th}).

\vspace{0.2cm}\noindent 
{\bf Surface cubications and  admissible immersions.}
The case of the surfaces is even simpler to understand. The immersion 
$\varphi_C(N_C)$ is obtained by  
drawing  arcs connecting the opposite sides for each square 
of the cubication $C$ and $N_C$ is a disjoint union of several circles.  
The immersions which arise from cubications are required to some mild 
restrictions. First the immersion is normal (or with normal crossings), since 
it has only transversal double points. All immersions encountered below 
will be normal crossings immersions.    
Since we can travel from one square of $C$ to any other 
square of $C$ by paths crossing the edges of $C$ it follows that 
the image of the immersion $\varphi_C(N_C)$ should be connected. 
On the other hand,  by cutting the surface $\Sigma$ along the arcs 
of $\varphi_C(N_C)$ we get a number of polygonal disks. An immersion 
having these two properties was called {\em admissible} in \cite{F}. 
Further we have a converse for the construction given above. 
If $j$ is an  admissible immersion of circles in the surface 
$\Sigma$ then $j$ is  $\varphi_C$ for some cubication $C$ of 
$\Sigma$. The abstract complex $C$ is the dual of the partition of 
$\Sigma$ into polygonal disks by means of the arcs of $j$. Since $j$ can have 
at most double points it follows that $C$ is made of squares.

\vspace{0.2cm}\noindent 
{\bf  Cubical flips on surfaces.} There are four different flips  (and their 
inverses) on  a surface, that we denoted by $b_1,b_2,b_3$ and $b_{3,1}$. 
They are pictured below. 

\begin{center}
\includegraphics{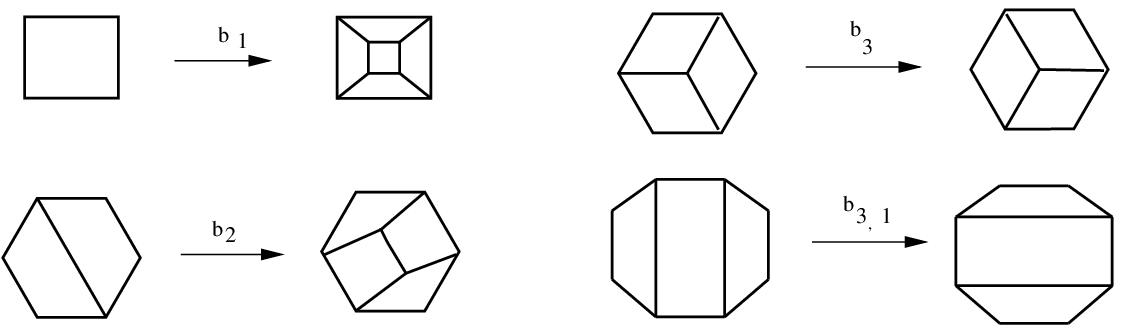}
\end{center}

\vspace{0.2cm}\noindent 
In particular, we have the flips denoted by the same letters that 
act on immersions of curves on surfaces. These transformations are {\em local moves}
in the sense that they change just a small part of the immersion 
that lives in a disk, leaving the immersion unchanged outside this 
disk. Specifically, here are the flip actions.

\begin{center}
\includegraphics{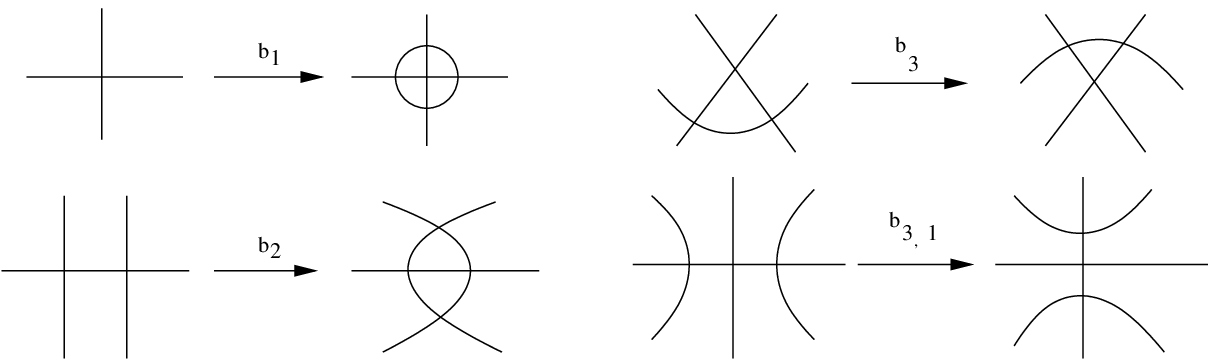}
\end{center}

\vspace{0.2cm}\noindent 
{\bf The invariant of cubications.} We associate to each  proper 
immersion 
$\alpha:L^1\to \Sigma$ of a disjoint union of circles and intervals $L^1$, 
two independent invariants, as follows. The image $\alpha(L^1)\subset \Sigma$
is a union of curves on $\Sigma$ and it can be viewed as a singular 
1-cycle of $\Sigma$. Notice that circles lay in the interior of $\Sigma$ 
while intervals are properly immersed and thus their endpoints lay on the 
boundary.  
We set then 
\[ j_1(\alpha)=[\alpha(L^1)]\in H_1(\Sigma,\partial \Sigma;\Z/2\Z)\]
Further we denote by $j_2(\alpha)\in \Z/2\Z$ the number of double points 
of $\alpha(L^1)$ mod two, and eventually 
\[ j_*(\alpha)=(j_1(\alpha),j_2(\alpha))\in H_1(\Sigma,\partial \Sigma;\Z/2\Z)\oplus \Z/2\Z\]
Further we are able to define the invariant associated to cubications 
by means of the formula: 
\[ j(C)=j_*(\varphi_C)\in  H_1(\Sigma,\partial \Sigma; \Z/2\Z)\oplus \Z/2\Z\]
Remark that we don't need to know that $j$ factors through the cobordism 
group $N(\Sigma)$ in order to define the invariant. 
Observe also that the boundary cubication is the disjoint union  
of polygons corresponding to the boundary circles and thus the numbers of 
edges determines completely their combinatorial type. 
The main theorem above is a consequence of the following 
more precise statement: 

\begin{theorem}\label{equiv}
Two marked cubications $C_0$ and $C_1$ of the compact surface 
$\Sigma$ are flip equivalent if and only if $j(C_0)=j(C_1)$ and 
their boundaries agree. 
\end{theorem}

\begin{remark} 
In the case of closed 
orientable surfaces our result is a consequence of the Nakamoto-Ota theorem 
(\cite{NO}). In fact, we will prove in section 5 that the diagonal  
transformations introduced by Nakamoto can be written as products of 
cubical flips. However, their method could not be used to cover the case where 
the surface is non-orientable or has boundary. 
Remark however that a weaker result holds true for  
diagonal transformations  on arbitrary closed surfaces 
(see \cite{Na1,Na2}), in which one  
replaced marked cubications up to isotopy by marked cubications up 
to homeomorphism. 
\end{remark}
\begin{remark}
Our methods are not combinatorial, 
as was the case of the sphere (see \cite{Th,F,BEE}),  
since one uses  in an essential 
manner the identification of $H_1(\Sigma,\partial \Sigma;\Z/2\Z)\oplus \Z/2\Z$ with 
$N(\Sigma)$,  which is of topological nature. 
The main interest in 
developing the topological proof below is that one can give an 
unifying treatment of all surfaces and the hope that  
these arguments might be generalized to higher dimensions.  
However, it would be interesting to find a direct combinatorial 
proof  that provides an algorithm which gives  
explicitly a sequence of flips connecting two cubications. 
Such an algorithm can be obtained in the closed orientable case by using 
the diagonal slides.     
\end{remark}

\vspace{0.2cm}\noindent 
Consider two cubications $C_0$ and $C_1$ having the same invariants. 
The first step in the proof of theorem \ref{equiv} is to return to the language of immersions and show that: 
\begin{proposition}\label{cobord}
If $C_0$ and $C_1$ have the same boundary and 
 $j(C_0)=j(C_1)$ then $\varphi_{C_0}$ and $\varphi_{C_1}$ are cobordant 
immersions. 
\end{proposition}
\noindent The second step is to use the existence of a cobordism  in order 
to produce flips and prove that: 
\begin{proposition}\label{flip}
If $\varphi_0$ and $\varphi_1$ are admissible immersions which are 
cobordant then there exists a sequence of flips which connect them. 
\end{proposition}
\noindent These two propositions will end the proof of the theorem \ref{equiv}.

\section{Cobordant immersions are flip equivalent}
\subsection{Flips, saddle and $X$-transformations relating cobordant immersions}
{\bf Connecting the immersions by means of maps with higher singularities.}
This section is devoted to the proof of proposition \ref{flip}. 
Consider thus a proper immersion 
$\varphi:F\to \Sigma\times [0,1]$ of a surface $F$ which is a cobordism 
between  the immersions $\varphi_0$ and $\varphi_1$. 

\vspace{0.2cm}\noindent 
The image $\varphi(F)$ is an immersed surface having therefore  a set of 
finitely many triples points that we denote by $S_3(\varphi(F))$. 
The set $S_2(\varphi(F))$ of double points of $\varphi$ (at the target) 
form a one-dimensional manifold, whose closure contains the 
triple points.  We have then a stratification of $\varphi(F)$ by 
manifolds 
\[ \varphi(F)=R(\varphi(F)) \cup S_2(\varphi(F)) \cup S_3(\varphi(F))\]
where $R(\varphi(F))$ is the set of non-singular (or regular) points.

\vspace{0.2cm}\noindent 
Our aim is to analyze the critical points of the restriction of the 
height function to $\varphi(F)$, by taking into account the singularities 
of $\varphi(F)$. In order to define critical points properly we need 
more terminology. Note that $R(\varphi(F))$ and $S_2(\varphi(F))$ 
are subsets of $\varphi(F)$ which might cause some troubles because 
critical points on the closure of a stratum might belong to another stratum.

\vspace{0.2cm}\noindent 
Any point $p\in \varphi(F)$ has an open neighborhood that is diffeomorphic 
to one coordinate plane, the union of two coordinate planes 
or the union of the three coordinate planes in $\R^3$, depending on 
whether $p\in R(\varphi(F)), p\in S_2(\varphi(F))$ or $p\in S_3(\varphi(F))$.  
The images of coordinate planes by this diffeomorphism are called the 
{\em leaves} of $\varphi(F)$ around $p$. Actually the leaves are well-defined 
only in a small neighborhood.  
A point $p\in \varphi(F)$ will be called {\em critical} for $h$ if 
$p$ is critical  either for the restriction of $h$ to some leaf 
containing $p$, or for $h|_{S_2(\varphi(F))}$, or else 
$p\in S_3(\varphi(F))$. Moreover, by using a small perturbation 
of $\varphi$ that is identity on the boundary, we can assume that 
the restriction of $h$ to the leaves is also a Morse function. 

\vspace{0.2cm}\noindent 
A consequence of the Morse theory is the following. If the 
interval $[t_1,t_2]$ does not contain any critical value for $h$ 
then there exists a diffeomorphism  of $\Sigma\times \{t_1\}$ 
into $\Sigma\times \{t_2\}$ 
that sends 
$\varphi(F)\cap h^{-1}(t_1)$ on 
$\varphi(F)\cap h^{-1}(t_2)$. 
Thus changes in the topology of the slice 
$\varphi(F)\cap h^{-1}(t)$ arise only at critical $t$. 

\vspace{0.2cm}\noindent 
A critical point $p$ will be said to be {\em unstable} if  
there are at least two leaves around $p$ and the restriction of 
$h$ to some leaf has a critical point at $p$. The other critical points 
are called {\em stable}.

\begin{lemma}
One can perturb slightly $\varphi$ by using an arbitrary 
small isotopy that is identity on the boundary such that 
$h$ has only stable critical points. 
\end{lemma}
\begin{proof}
Let $p$ be an unstable critical point and $\alpha$ be a leaf 
that is critical for the restriction of $h$ at $p$. 
This is equivalent to the fact that the gradient of $h$ 
(which is nowhere zero since $h$ is regular) 
is orthogonal to the tangent plane at 
$\alpha$. Since the immersion is normal crossings 
any other leaf $\beta$  around $p$ should be transverse to $\alpha$ and 
thus the restriction of $h$ to that leaf is non-critical at $p$. 

\vspace{0.2cm}\noindent 
Use now a small perturbation of the extra leaf around $p$ by an isotopy 
that moves the intersection arc $\alpha\cap \beta$  
off $p$. The new intersection point is not 
anymore critical for the restriction of $h$ at $\alpha$.  
\end{proof}
 
\vspace{0.2cm}\noindent  
There are the following situations when the values are not regular:
\begin{enumerate}
\item the slice $\Sigma\times\{t\}$ passes thru a triple point $p$ 
and the restriction of $h$ to all leaves is regular. 
\item the slice $\Sigma\times\{t\}$ contains a critical point $p$ 
of the restriction  $h|_{S_2(\varphi(F))}$, 
to the double points locus. Let $\alpha$ and $\beta$ denote 
the two leaves around $p$. The restriction of 
$h$ to the two leaves is regular. 
\item If the slice  $\Sigma\times\{t\}$ contains a critical point
of the height restriction $h|_{R(\varphi(F))}$, to the regular locus. 
\end{enumerate}

\vspace{0.2cm}\noindent 
{\bf Passing a triple point.} We can assume by general position arguments 
that the (finitely many) triple points have  distinct heights and thus 
the slice $\Sigma\times\{t_c\}$ contains precisely one triple point $p$ and 
no other critical point.

\begin{proposition}
When crossing a critical value $t_c$  corresponding to a stable triple point  
the image of the  sliced immersion $\varphi(F)\cap \Sigma\times \{t\}$ 
changes according to a flip $b_3$. 
\end{proposition}
\begin{proof}
Consider a coordinates chart $(V, \nu)$ on $\Sigma \times [0,1]$  containing 
the point $p$. We assume that $V$ is diffeomorphic to $\R^3$ and the 
diffeomorphism $\nu$ sends $\varphi(F)$ into $\Pi$, where  $\Pi$ denotes 
the union of the three coordinates planes in $\R^3$. 
Denote by $H:\R^3\to \R$ the function $h$ expressed in these coordinates. 
Notice that $H$ is regular and one can assume that $H(0)=0$.  
The level hypersurface ${\mathcal H}=H^{-1}(0)$ corresponds to a 
neighborhood of $p$ into the critical slice. Moreover, one knows that 
$H$ has no critical points when restricted to the leaves around $p$ 
(i.e. the coordinate planes) which amounts to say that 
the gradient ${\rm grad}_0 H$ is not orthogonal to any coordinate 
plane. 
Let $T_0{\mathcal H}$ be the tangent space at ${\mathcal H}$ at the origin.
There is only one generic position of a plane through $O$ with respect to 
$\Pi$. Translating the plane along the direction of ${\rm grad}_0 H$ yields 
the desired bifurcation.  
\end{proof}

\vspace{0.2cm}\noindent 
{\bf Critical points on the double points locus.}
The double locus $S_2(\varphi(F))$  is a 1-manifold and each of its  points  
has a coordinate chart in which 
$\varphi(F)$  is sent into   
the union of  the first two 
coordinates planes in $\R^3$. The critical points correspond to 
local extrema of $h$ when restricted to the double line.  

\begin{proposition}
When crossing a critical value $t_c$  corresponding to a stable critical 
point of the restriction  of $h$ to  $S_2(\varphi(F))$   
the image of the  sliced immersion $\varphi(F)\cap \Sigma\times \{t\}$ 
is transformed according to the following picture: 
\begin{center}
\includegraphics{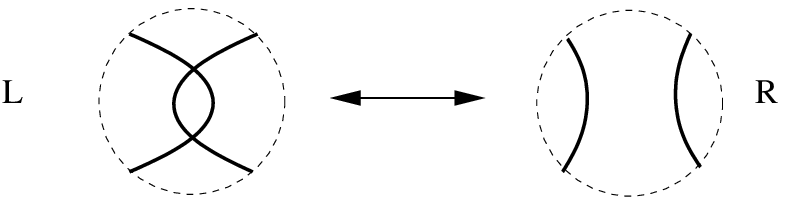}
\end{center}
This means that we replace the left diagram $L$ by the right diagram $R$, 
leaving the part of the immersion unchanged outside the small ball 
figured above. 
\end{proposition}
\begin{proof}
The proof from above applies with minor modifications. 
However here the singularity 
consists of  a line of double points and thus the intersection 
of the two planes with the boundary sphere of a small disk centered at 
the critical point is the union of two great 
circles instead of three. Since there is an extremum  on the double line 
there is one half of the sphere that contains no double points. Thus the 
only transformation  possible is that from the picture.  
\end{proof}

\vspace{0.2cm}\noindent 
The local transformation  of an immersion (and its inverse) 
occurring in the proposition above is called an {\em $X$-transformation}.

\vspace{0.2cm}\noindent 
{\bf Critical points on the regular strata.}
The height function $h|_{R(\varphi(F)}$ is a Morse function on a surface 
and thus it has critical points of three types: maximum, minimum and 
saddle (index one). 

\vspace{0.2cm}\noindent 
Passing thru a minimum point amounts to adjoin a small embedded sphere 
to the immersion, while a maximum point contributes with deleting a small 
sphere. We will speak about creation/annihilation of a circle.

\vspace{0.2cm}\noindent 
When passing thru a saddle point the  
sliced immersion $\varphi(F)\cap \Sigma\times \{t\}$ is subject to the  
following familiar change:

\begin{center}
\includegraphics{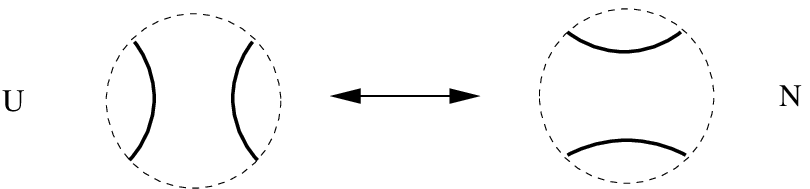}
\end{center}
\noindent 
This transformation will be called {\em saddle transformation}.

\subsection{Stabilizations}
In order to establish the claim  of proposition \ref{flip} it would suffice to show that 
any saddle move, $X$-move or creation/annihilation can be actually 
realized by means of flips. As stated, this cannot be true, but a slight 
modification of this statement will hold true. One  reason is that 
creating a new circle destroys the connectivity of the graph, while  
flips preserve it. Recall however that our immersions were supposed to be 
admissible. In particular, in dimension 2 one has a connected graph 
whose complementary is a union of open disks on the surface. 

\vspace{0.2cm}\noindent 
We will introduce now another operation on  immersions (or, equivalently 
on graphs on surfaces) that will be called {\em stabilization}. 
Let $\varphi:\sqcup S^1\to \Sigma$ be an immersion of circles and 
$\lambda$ a simple arc in $\Sigma$ that joins two 
points of the image of $\varphi$ and is transversal to it. 
Let $u(\lambda)$ be the embedding 
of a circle into $\Sigma$ as the boundary 
of a small regular neighborhood of $\lambda$. Then the union of 
$\varphi$ and $u(\lambda)$ define an immersion that we denote by 
$\varphi*_{\lambda}u$. One allows also the  particular case 
when $\lambda$ is trivial and the two points coincide, denoted 
$\varphi *u$. We will also denote by $\varphi \sqcup u$ the disjoint 
union of $\varphi$ with a small trivially embedded circle.  
 
\vspace{0.2cm}\noindent 
The main result of this subsection is: 

\begin{proposition}\label{stab} 
 If $\varphi$ is admissible then $\varphi*_{\lambda}u$ is flip equivalent to 
$\varphi$. 
\end{proposition}
\begin{proof}
Assume that $\lambda$ intersects the image of  
$\varphi$ at the points $p_1,p_2,...,p_k$, where $k\geq 2$. 
We use induction on $k$. Let $\lambda'$ be the sub-arc 
of $\lambda$ joining the last two points  $p_{k-1}$ and $p_k$ and 
$\lambda"$ its complementary arc. Then $u(\lambda)$ is the 
connect sum of $u(\lambda')$ and $u(\lambda")$. Let $z$ be the arc which 
is common to $u(\lambda')$ and $u(\lambda)$. 
We will use flips whose action is trivial outside the arc $z$.

\vspace{0.2cm}\noindent 
The circle $u(\lambda')$ intersects 
two arcs $a$ and $b$ of  the image of $\varphi$ in two points 
that determine an arc $q$ of 
$u(\lambda)$ which is free of other intersection points.
The arcs $a$ and $b$ are connected by some chain of arcs of $\varphi$, namely  
$a_1,a_2,...,a_m$, since $\varphi$ is admissible. 
Moreover, the arcs $a_i$, $a, b$ and $q$  
are bounding a face of the complementary of the immersion. 
Since $\varphi$ was supposed admissible the  simple arc $q$ 
subdivides a topological disk into two pieces which should be 
again disks. Thus the polygon determined by  $a_i$, $a, b$ and $q$ 
is  a topological disk $Q$. 

\vspace{0.2cm}\noindent 
Now, we can slide the intersection part between $z\cap a$ 
along the path $a_1,a_2,...,a_m,b$ until $z$ will intersect 
both $a_m$ and $b$ as below: 

\begin{center}
\includegraphics{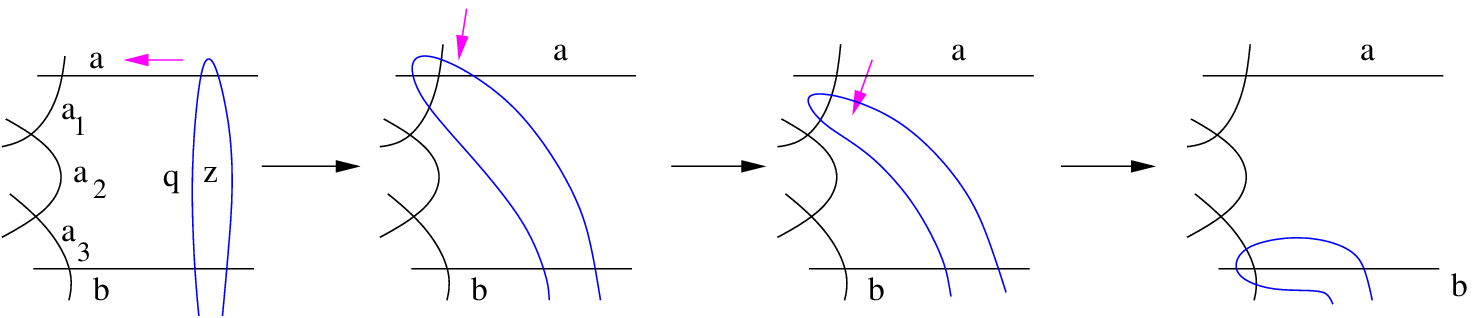}
\end{center}

\vspace{0.2cm}\noindent 
This can be realized by flips that acts non-trivially only on the 
arc $z$ and are identity outside. 

\vspace{0.2cm}\noindent  
If $k=2$ then $u(\lambda)$ is obtained from $z$ by joining 
its extremities. Then the inverse flip move 
$b_1$ destroys the circle $u(\lambda)$.

\vspace{0.2cm}\noindent  
If $k >2$ then use the inverse move $b_2$. The result of this sequence 
of flips is that $u(\lambda)$ is transformed 
into $u(\lambda")$ and $\lambda"$ has $k-1$ intersection points with the 
image of $\varphi$. The claim follows then by induction.  
\end{proof}

\subsection{Saddle transformations}
Consider an immersion $\psi$ containing a small ball where it coincides 
with the diagram $U$, the left hand side in the picture of the saddle move.
We denote by $S(\psi)$ the immersion obtained  by a saddle move from $\psi$, 
namely by excising $U$ and gluing back $N$. 
\begin{lemma}\label{sad}
Assume that the immersion $\psi$  is admissible. 
Then $\psi$ is flip equivalent to  some stabilization $S(\psi)*_{\lambda}u$
of $S(\psi)$. 
\end{lemma}
\begin{remark}
The fact that in general $\psi$ is equivalent only with the stabilization 
$S(\psi)$ is not unexpected. In fact saddle move 
might destroy the connectivity of the image, and for 
that reason one should add a new circle $u(\lambda)$ in order 
to restore it.  
\end{remark}
\begin{proof}
It is essential that $\psi$ is admissible and thus the image is connected. 
In particular, there are several arcs of $\psi$ which join the 
left arc of $U$ to the right arc of $U$. These intermediary arcs are 
outside the small ball. However, we can choose  a specific family of arcs, 
namely those bounding a face (i.e. a topological disk) 
in the complement of  the image of $\psi$. 

\vspace{0.2cm}\noindent 
We use inductively the move $b_2$ in order to push and slide 
the left arc of $U$ across each intermediary arc, 
as it is shown below for the case of the first arc. 

\begin{center}
\includegraphics{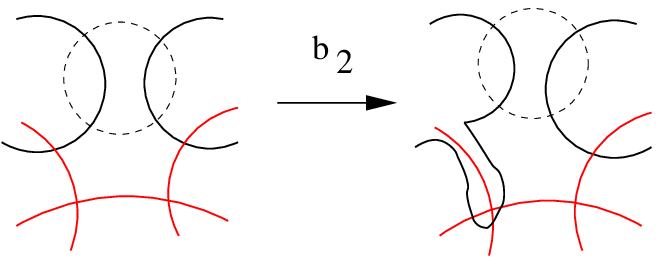}
\end{center}

\noindent We do that until the new position of the left arc 
intersects the last intermediary arc. Use $b_1$ to create a small 
circle centered at the intersection point between the last intermediary arc 
and the right arc of $U$. The next step is to use the move $b_{3,1}$ 
once. Its  support is centered at  the intersection between 
the last intermediary arc and the small circle created at the previous step. 
We obtain the configuration below:

\begin{center}
\includegraphics{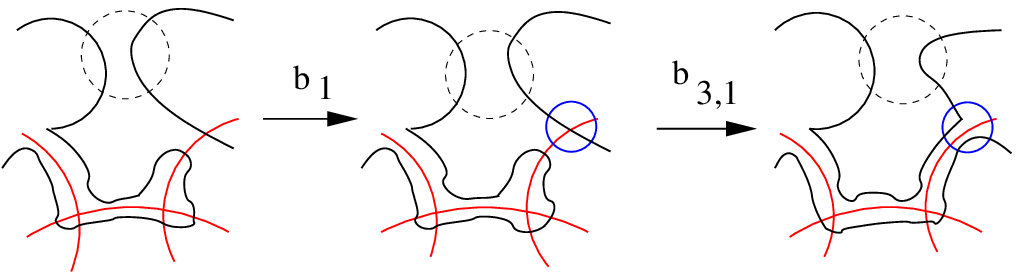}
\end{center}

\noindent 
We realized half of the saddle move, namely the upper arc from $N$. 
However, the bottom arc intersects all intermediary arcs used above. 
We will further use the inverse moves $b_2$ in order to slide the 
bottom arc over the intermediary arcs, this time from the 
right side back to the left. We first use the intersections with the 
additional small circle and then go along the intermediary arcs 
as follows:

\begin{center}
\includegraphics{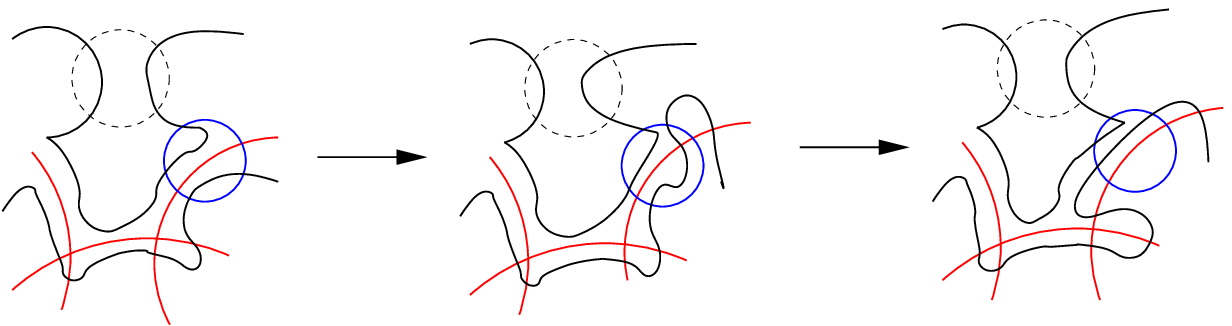}
\end{center}

\noindent 
At the end we transformed the bottom arc into another arc which 
go along the path of intermediary arcs but lays on its upper side. 
Further the small additional circle can be slid outside the 
last intermediary arc. The face determined by the intermediary arcs 
is topologically a disk, since the immersion is admissible. Thus 
both arcs (upper and lower) can be isotoped to the position that 
they have in the $N$ diagram.  We obtained  thus 
the following configuration:

\begin{center}
\includegraphics{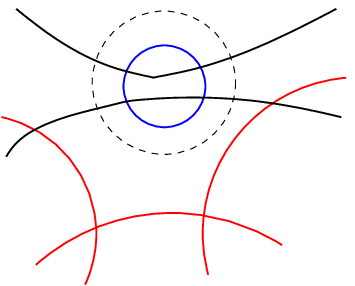}
\end{center}
 
\noindent This is the stabilization  $S(\psi)*_{\lambda}u$
where we used the arc $\lambda$ that connects the two 
arcs in the diagram $N$. 
\end{proof}

\subsection{$X$ transformations}
We denote again by $\psi$ an immersion containing  in a small 
ball the diagram $L$ and 
by $X(\psi)$ the result of the $X$ transformation. 
\begin{lemma}\label{X}
If $\psi$ is admissible then $\psi$ is flip equivalent to 
$X(\psi)*_{\lambda}u$. 
\end{lemma}
\begin{proof}
We use $b_1$ to create a small 
circle centered at one intersection point of the arcs in $L$. 
Using next the move $b_2$ we can slide the two arcs from $L$ far apart 
and obtain $X(\psi)*_{\lambda}u$.

\begin{center} 
\includegraphics{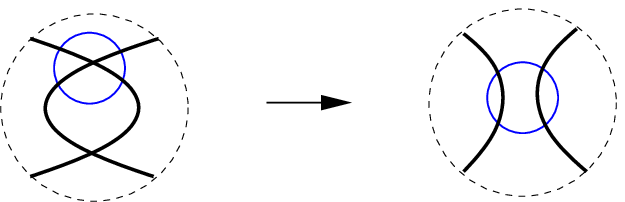}  
\end{center}
\end{proof}

\vspace{0.2cm}\noindent Consider now the immersion $\psi$ that contains 
the diagram $R$ so that we can apply the inverse move $X^{-1}$ to 
the immersion $\psi$. In this case we have 

\begin{lemma}
If $\psi$ is admissible then $X^{-1}(\psi)$ is flip equivalent to $\psi$. 
\end{lemma}
\begin{proof}
The figure above shows that $X^{-1}(\psi)$ is flip equivalent 
to $\psi*_{\lambda}u$. If $\psi$ is admissible then $\psi$ 
is equivalent to $\psi*_{\lambda}u$ by lemma  \ref{stab}. 
\end{proof}

\subsection{Proof of proposition \ref{flip}}
Consider now two  admissible immersions $\varphi$ and $\psi$ 
which are cobordant. According to the previous description of 
cobordisms there exists a sequence  of  moves 
$A_1,A_2,...,A_{n-1}$  transforming 
$\varphi$ into $\psi$ which are 
either flips, saddle transformations, $X$  moves or their inverses, 
or creation/annihilation moves. Let the sequence of 
immersions so obtained be 
$\varphi=\varphi_1,\varphi_2,\ldots, \varphi_n=\psi$. 

\vspace{0.2cm}\noindent
The problem we face is that some of the $\varphi_j$ 
might have disconnected image and are therefore not admissible. 
First, we can dispose of the creation/annihilation of circles, by using 
instead saddle moves, as below: 

\begin{center}
\includegraphics{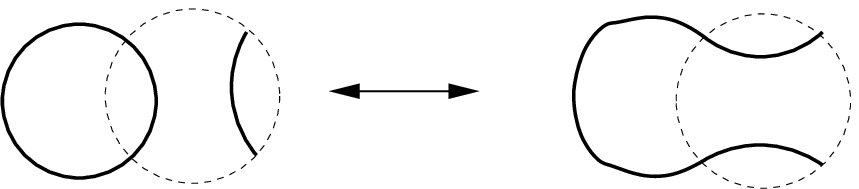}
\end{center}
 
\vspace{0.2cm}\noindent
The key ingredient is the existence of an admissible 
modification of the sequence above. Let us introduce
first more terminology. Assume that we have a sequence of moves $A_j$ that 
might be saddle, $X$-transformations of flips. 
The move $A_j$ replaces the part of the image  of some immersion 
$\varphi_j$ contained in a disk $D_j$ by a different graph, while 
keeping the complementary unchanged. 
We define now an {\em extended move} $\overline A_j$ that acts on 
some stabilization $\overline{\varphi_j}$ of $\varphi_j$ that uses 
the arcs $\lambda_i$, as follows.  
The extended move $\overline{A_j}$ is assumed to act exactly 
in the same way as $A_j$, namely it replaces the part of the 
image of $\varphi_j$ contained in the disk $D_j$ by the corresponding 
graph (as prescribed by $A_j$), while keeping untouched both the 
complementary of the ball and also the new circles $u(\lambda_i)$.
The extended moves $\overline A_j$ are not anymore 
usual saddles, $X$-moves etc unless the disk $D_j$ avoids 
the  arcs $\lambda_i$. Notice also that  
$\overline A_j\overline \varphi_j$ might be non-admissible.
\begin{lemma}
There exists a sequence of immersions  $\overline \varphi_j$ 
and extended moves $\overline A_j$ with the following properties: 
\begin{enumerate}
\item $\overline \varphi_j$ is a stabilization of $\varphi_j$, of the 
form 
$ \overline \varphi_j= {\varphi_j}* _{\lambda_{1}}u*_{\lambda_{2}}u
\cdots*_{\lambda_{j-1}}u$. 
Thus $\varphi_j$ is obtained by $\varphi_j$ by adding several circles 
of type $u(\lambda_i)$, where the arcs $\lambda_i$ might intersect 
each other. 
\item  The immersions $\overline \varphi_j$ are admissible and  
$\overline \varphi_j$ and $\overline\varphi_{j+1}$ are flip equivalent.
\end{enumerate}
\end{lemma} 
\begin{proof}
We use a recurrence on  the length of the sequence. 
By the definition of the extended moves we have:  
\[ \overline A_j \overline \varphi_j= 
 {\varphi_{j+1}}* _{\lambda_{1}}u*_{\lambda_{2}}u
\cdots*_{\lambda_{j-1}}u, {\rm and}~ 
\overline \varphi_{j+1}=({\overline A_j\overline \varphi_j})* _{\lambda_j}u\]
\vspace{0.2cm}\noindent 
We set  $\overline\varphi_1=\varphi_1$.
Then  $A_1\in\{ X,X^{-1},S,b_3\}$. We further define
\[\overline{\varphi_2}= \left\{\begin{array}{ll}
\varphi_2, & {\rm  if } ~  A_1\in \{X^{-1},b_3\} \\
\varphi_2 *_{\lambda_1} u, & {\rm otherwise}  
\end{array}\right.
\] 
where the arc $\lambda_1$ is that furnished by the lemmas 
\ref{sad} and \ref{X}. According to these lemmas 
$\overline{\varphi_2}$ is flip equivalent to $\overline{\varphi_1}$ and 
both are admissible. Moreover, $\overline{A_1}=A_1$. 

\vspace{0.2cm}\noindent Assume now that $\overline{\varphi_j}$ and 
$\overline{A_j}$ are defined for $j\leq k$. 
If $A_k\in \{b_3, X^{-1}\}$ then we set 
\[\overline{\varphi_{k+1}}=\overline{A_k}\overline{\varphi_k}\]
In the other two cases we have to analyze the picture inside the 
disk $D_k$.  We set first 
$\Lambda=D_k \cap (\lambda_1 \cup \lambda_2 \cup \cdots \cup \lambda_{k-1})$  
and call its components the special arcs.  

\begin{enumerate}
\item Assume that $A_k=X$. We can take  for $D_k$ a very tiny  
regular neighborhood of the two arcs $a$ and $b$ of $\varphi$  
joining the double points in the diagram $U$ 
and thus all arcs $\lambda_j$ entering $D_k$  should intersect 
these arcs. Another useful observation is that any sub-arc of a special 
arc that has two consecutive  intersection points with $a$ can be 
moved off $D_k$ by means of $b_2$ moves. Thus either there are no special arcs 
within $D_k$ or else there exists arcs that cross both arcs $a$ and $b$.

\begin{enumerate}
\item If $\Lambda=\emptyset$ then we set $\overline{\varphi_{k+1}}=
X(\overline{\varphi_k})*_{\lambda_k}u$ where $\lambda_k$ 
is the arc given by lemma \ref{X}. That lemma shows that 
$\overline{\varphi_{k+1}}$ is flip equivalent to 
$\overline{\varphi_{k}}$. 
\item
Otherwise we define  $\overline{\varphi_{k+1}}=
\overline{A_k}(\overline{\varphi_k})$, and it suffices to show the 
flip equivalence with $\overline{\varphi_k}$. From above one knows that 
$ a\cup b\cup \Lambda$ is  
connected. We can use the moves $b_3$ in order to move the arc $a$ 
across the vertices of the diagram $\Lambda$ and the moves 
$b_2$ in order to move the arc $a$ along special arcs without vertices on it
i.e. arcs which are connected components of $\Lambda$.
\begin{center}
\includegraphics{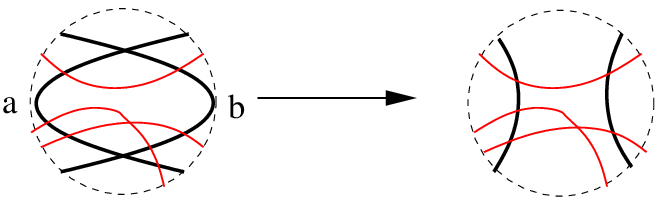}
\end{center}
\end{enumerate}
\item Consider now that $A_k=S$. Let $a$ and $b$ denote the arcs 
from the diagram $U$ and $a',b'$ the arcs from the diagram $N$. 
\begin{enumerate}
\item If $\Lambda=\emptyset$ then we set $\overline{\varphi_{k+1}}=
S(\overline{\varphi_k})*_{\lambda_k}u$ where $\lambda_k$ 
is the arc given by lemma \ref{sad}. That lemma shows that 
$\overline{\varphi_{k+1}}$ is flip equivalent to 
$\overline{\varphi_{k}}$. 
\item The move $S$ takes place in a very tiny neighborhood of an
arc that joins two points, one from $a$ and the other from $b$. 
We call it the {\em core} and it corresponds to the 1-handle 
to be added to the immersion by the saddle move. 
Thus we can discard (by using an isotopy) all special arcs except 
those that intersect the core. Moreover, one can suppose that there 
is at least one such special arc.  We define then 
$\overline{\varphi_{k+1}}=
\overline{A_k}(\overline{\varphi_k})$, and it suffices to show the 
flip equivalence with $\overline{\varphi_k}$.
Lemma \ref{sad} shows that $\overline{\varphi_k}$ 
is flip equivalent to  $\overline{A_k}(\overline{\varphi_k})*_{\lambda_k}u$, 
where $\lambda_k$ is the arc that joins $a'$ to $b'$. One has to notice that 
$\Lambda\cup a'\cup b'$ is connected since the special arcs intersect 
the core. Thus a suitable application of the moves $b_3$, slidings across 
the special arcs and an inverse $b_1$ move will get rid of the 
extra circle $u(\lambda_k)$. 
\begin{center}
\includegraphics{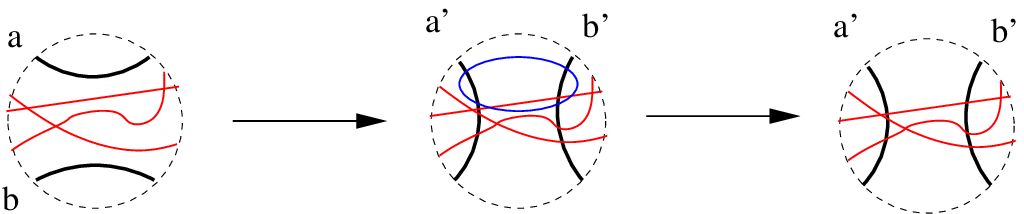}
\end{center}
\end{enumerate}
\end{enumerate}
\end{proof}

\vspace{0.2cm}\noindent 
Now the lemma ends the proof of the 
proposition \ref{flip} immediately. In fact one knows  that 
the last immersion $\varphi_n$ is admissible, by hypothesis. 
Then lemma \ref{stab} stated that 
$\overline{\varphi_n}=\varphi_n*_{\lambda_1}u*_{\lambda_2}u*\cdots 
*_{\lambda_{n-1}}u$ is flip equivalent to $\varphi_n$. On the other hand 
$\varphi_1$ is equivalent to $\overline{\varphi_j}$ for all $j$, 
and we are done.

\section{The cobordism group of immersions in a surface}
Although the computation of $N(\Sigma)$ for closed surfaces 
is folklore we didn't find a reference addressing precisely this issue.  
Related results of similar nature are recorded in \cite{Ca,CS}.  
This also follows from our computations of the cobordism groups of 
codimension one immersions for  manifolds of  small 
dimensions, from \cite{FG}.  
The proof that we give below is elementary  in that it 
uses only basic methods and results from three dimensional topology. 
The proposition \ref{cobord} can be reformulated as follows: 
\begin{proposition}
Two  one-dimensional immersions in a compact surface are cobordant 
if and only if they have the same  boundary and their $j_*$ invariants agree.  
In particular the map $j_*$ factors to an isomorphism 
$j_*:N(\Sigma)\to H_1(\Sigma,\partial \Sigma;\Z/2\Z)\oplus \Z/2\Z$. 
\end{proposition}
\begin{proof}
Assume that $\varphi_i:L_i^1\to \Sigma$, $i=0,1$ are proper immersions 
of two disjoint unions of circles and intervals, $L_0^1$ and $L_1^1$, 
having the same boundary and 
invariants. Thus  the number of intervals within 
$L_i^1$ is the same for both values of $i$, and their images  
by the respective immersions intersect $\partial \Sigma$ 
in the same set of points. 

\vspace{0.2cm}\noindent 
Since $j_1(\varphi_0)=j_1(\varphi_1)$ there exists a 
singular 2-cycle of  the pair 
$(\Sigma,\partial \Sigma)$ with  boundary 
$\varphi_1(L_0^1)-\varphi_2(L_1^1)$. 
This means that there exists a surface with corners $F$ having boundary 
$\partial F=L_0^1 - L_1^1$ and a (singular) map 
$f:F\to \Sigma$ such that 
$f |_{\partial F}=\varphi_0 \sqcup \varphi_1$.
The corners set is the union of boundary points of intervals in 
$\varphi_i(L_i^1)$.

\vspace{0.2cm}\noindent 
Let $h:F\to [0,1]$ be a proper smooth function having $h^{-1}(0)=L_0^1$ and 
$h^{-1}(1)=L_1^1$.  We can lift $f$ to a  map 
$\phi:F\to \Sigma \times [0,1]$ by means of the formula 
\[ \phi(x)=(f(x), h(x))\]
Then $\phi$ is a proper map i.e. it sends the boundary into the 
boundary. 

\vspace{0.2cm}\noindent 
Let us consider first the case when the surface $\Sigma$ has no boundary 
and thus there are no intervals among the components of $L_i^1$. 
Thus $F$ is a surface with boundary.

\begin{lemma}
Let $\phi:F\to M^3$ be a proper map of a surface 
into the connected 3-manifold $M^3$, whose restriction 
to the boundary $\partial \phi:\partial F\to \partial M^3$ is an immersion.
Then  there exists an immersion $\phi':F'\to M^3$ of a possibly different 
surface $F'$ such that $\partial F=\partial F'$ and  
$\partial \phi'=\partial \phi$ if and only if 
$j_2(\partial \phi)=0$. 
\end{lemma}
\begin{proof}
The proof is a variation of that given by Hass and Hughes 
(\cite{HH}) in the case when 
$F$ is closed and thus the invariant $j_2$ trivially vanishes. 
The main arguments below appeared also in Whitehead's paper  
(\cite{Wh2}, proof of Theorem 4.1). 

\vspace{0.2cm}\noindent 
We use a classical result of Whitney which states that any smooth map $\phi:F\to M$ 
which is transverse to the boundary 
is homotopic rel boundary to a proper  general position map having only 
{\em simple branch points}. A  simple branch point 
in the image is a point having a 
neighborhood  homeomorphic to Whitney's umbrella, 
namely the cone over the bouquet of two 
circles (usually known as the figure eight), as pictured below: 

\begin{center}
\includegraphics{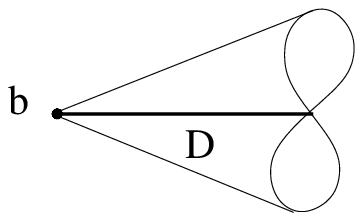}
\end{center}

\vspace{0.2cm}\noindent
According to Thom this is the only generic local singularity of 
smooth germs of maps  $(\R^2,0)\to (\R^3,0)$. 

\vspace{0.2cm}\noindent
In particular,  from any branch point $b$ of $\phi(F)$ emerges 
a line of double points, denoted by  $D$ in the figure above. 
Thus the map is now an immersion everywhere except at the (finitely many)
branch points. 

\vspace{0.2cm}\noindent 
The singularities of the new map (that we denote by the same 
letter $\phi)$ consist of (see also \cite{H}, p.41-42): 
\begin{itemize}
\item finitely many triple points, each one having three points in the 
preimage;
\item finitely many branch points, each one having one point in the preimage;
\item curves of double points, which are points where the surface 
$\phi(F)$ has normal crossing self-intersections. These curves are of 
two types:  
\begin{itemize}
\item type I curves that are either closed simple loops or 
joining two boundary double points of $\partial \phi$. These 
curves intersect transversally at triple points;   
\item type II curves that are either segments which join 
two distinct branch points or else a branch point to some 
boundary double point of $\partial \phi$. 
\end{itemize}
\end{itemize}

\vspace{0.2cm}\noindent 
If the map $\phi$ is an immersion then there are no branch points and thus 
the double points of the restriction $\partial \phi$ are paired by means 
of the type I curves, and thus there is an even number of such double points. 
This establishes the necessity of our condition. 

\vspace{0.2cm}\noindent
Conversely, if there is an even number of such double points, then 
we have a number of branch points paired together by means of 
type II curves and also an even  number of pairs $(b_j,p_j)$ where 
$b_j$ is a branch points and $p_j$ is a boundary double point.

\vspace{0.2cm}\noindent 
Let $b_1,b_2$ being two branch points connected by a curve of type II. 
There exists a standard 
way to modify $\phi$ by means of a homotopy so that the 
branch points are pushed one 
end towards the other (see also \cite{HH}):

\begin{center}
\includegraphics{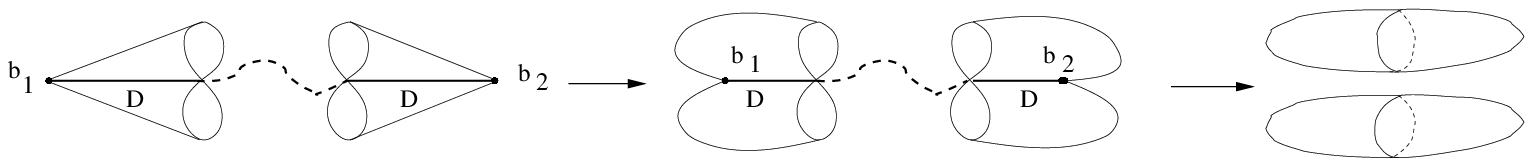}
\end{center}

\vspace{0.2cm}\noindent 
This is equivalent to cutting the surface $\varphi(F)$ along 
the connecting double curve. The domain of the map will change 
accordingly  but we can further add a small tube  between the components 
separated by the cut, in order to recover the same surface $F$.

\vspace{0.2cm}\noindent 
The other case is when we have a couple $(b_1,p_1)$ and $(b_2,p_2)$ 
of branch points paired with boundary  double points. Choose now 
a segment embedded in $M$ that joins $b_1$ to $b_2$ 
and avoids the triple points and other branch points. 
We will push again the branch points towards one another along this segment 
until they collide and then disappear, as in the figure below: 

\begin{center}
\includegraphics{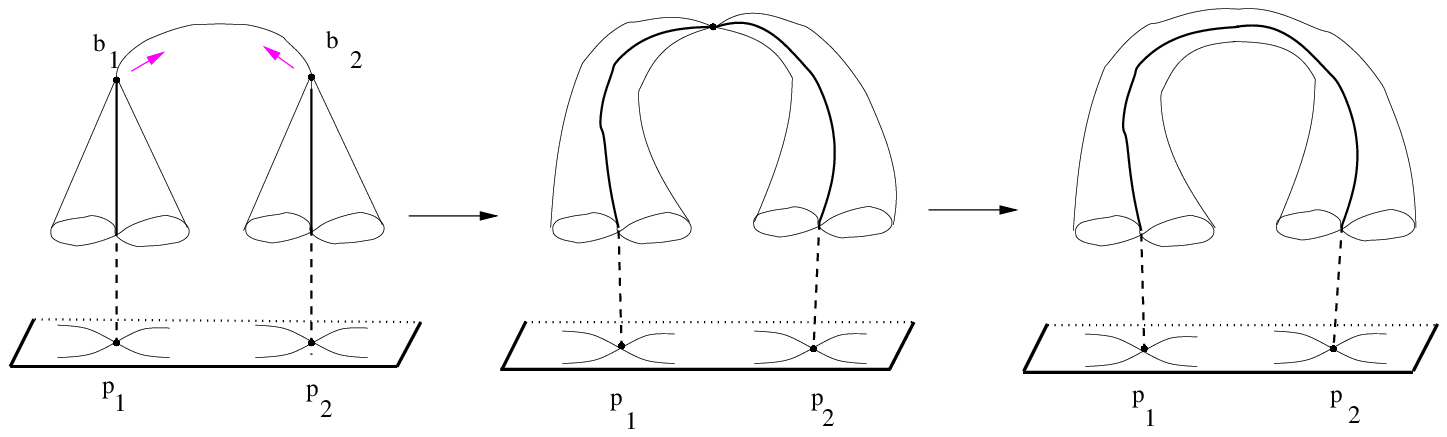}
\end{center}

\vspace{0.2cm}\noindent 
In particular $\phi$ was changed into an immersion $\phi'$ of a surface 
$F'$, which is obtained from $F$  by adding a 1-handle.

\vspace{0.2cm}\noindent 
Eventually we can slightly perturb the immersion in order to become 
normal crossings immersions, since the later are generic. 
\end{proof}
\noindent The lemma settles our claim, because we know that 
$j_2(\phi)=j_2(\varphi_0)-j_2(\varphi_1)=0$  and thus there exists 
an immersion $\Phi:F\to \Sigma\times[0,1]$ extending the immersions 
$\varphi_0\sqcup \varphi_1$ on the boundary.  
\end{proof}

\begin{remark}
It is considerably more difficult to get control on the genus of $F'$ 
out of the data $\phi,F,M$. The typical result is the celebrated 
Dehn lemma (\cite{H}). 
\end{remark}

\vspace{0.2cm}\noindent Consider now the general case of surfaces
$\Sigma$ with nonempty boundary, which is more difficult only at the 
terminology level.  The main difference is that now $\Sigma\times [0,1]$ 
is a manifold with corners. Then 
the map $\varphi:F\to \Sigma\times [0,1]$ is a proper 
map between manifolds with corners that respects the horizontal 
part of the boundary and it is a product on the vertical boundary. 
 
\begin{lemma}
Let $\phi:F\to M^3$ be a proper map of a surface  with corners 
into the connected 3-manifold $M^3$, whose restriction 
to the boundary $\partial \phi:\partial F\to \partial M^3$ is an immersion 
which respects the horizontal 
part of the boundary and it is a product on the vertical boundary. 
Then  there exists an immersion  between manifolds with corners 
$\phi':F'\to M^3$ of a possibly different 
surface with corners $F'$ such that $\partial F=\partial F'$ and  
$\partial \phi'=\partial \phi$ if and only if 
$j_2(\partial \phi)=0$. 
\end{lemma}
\begin{proof}
The proof from above applies with only minor modifications. 
In fact the lines of double points of $\phi$ in $M^3$ are disjoint from the 
vertical boundary $\partial_VM$ because the immersion is proper 
and its restriction $\partial \phi$  is a product on the vertical boundary.  
\end{proof}

\section{Diagonal transformations}
Nakamoto and Ota (see \cite{Na1,Na2,NO}) considered the problem 
of moves acting on the 
cubications of surfaces but they used a different family of 
transformations that they called diagonal transformations. 
Specifically, these are the diagonal slide $R$ (or rotation of order three) 
and  the diagonal rotation $D$ (or order two rotation) drawn below: 

\begin{center}
\includegraphics{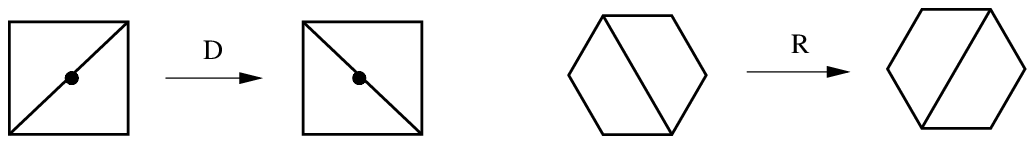}
\end{center}

\begin{proposition}
The diagonal transformations are products of cubical flips. 
\end{proposition}
\begin{proof}
There is an immediate corollary of the main result for surfaces 
with boundary since the immersions associated to these 
cubications of the disk have the same number of double points 
and boundary points. However, there exists an elementary 
pictorial  proof, in which we decompose $D$ into flips as below:

\begin{center}
\includegraphics{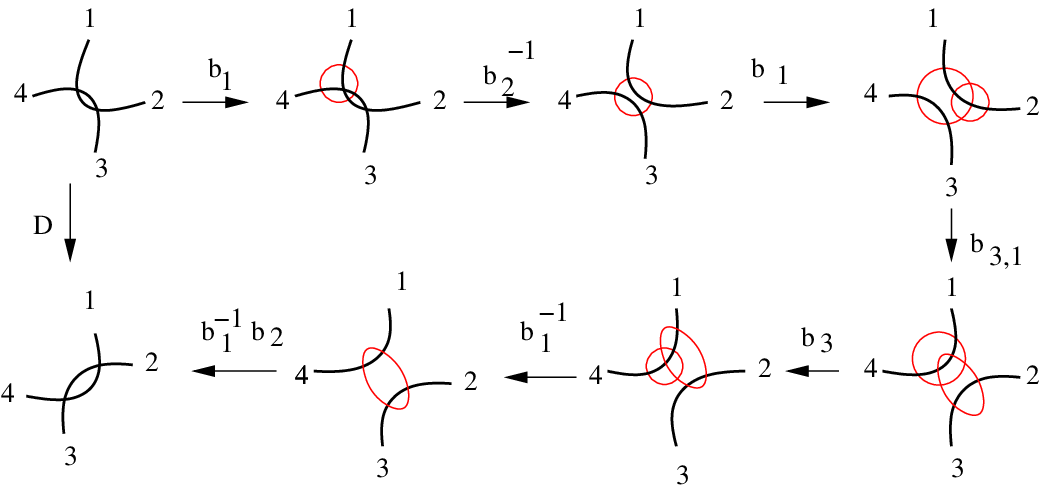}
\end{center}

\noindent and further the order three rotation move $R$: 
\begin{center}
\includegraphics{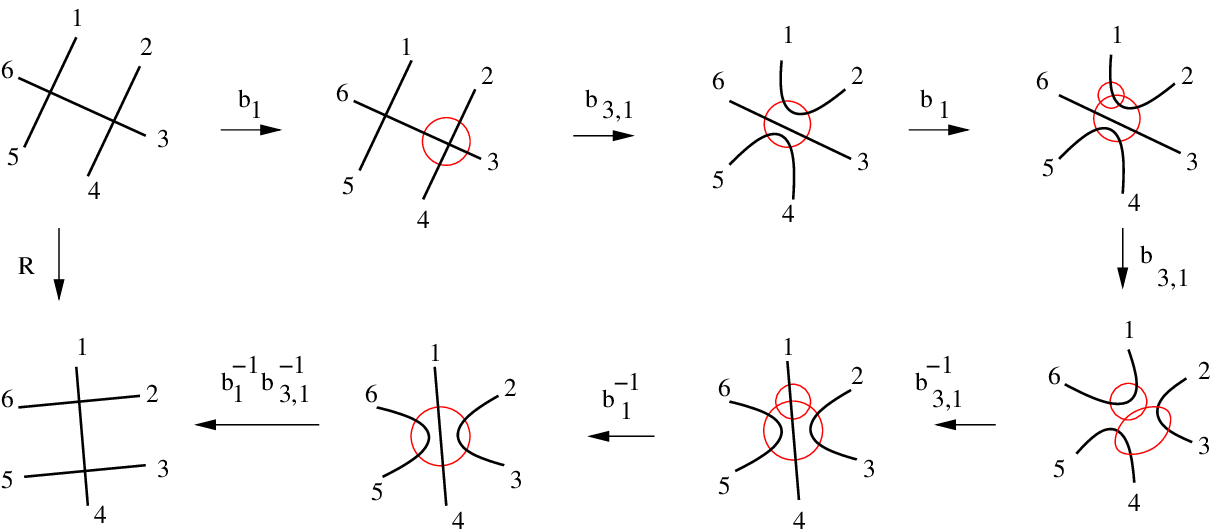}
\end{center}

\end{proof}

\noindent 
Nakamoto and Ota proved in \cite{NO} that  
there exists a natural number  $n_g$  such that any two cubications 
of a closed orientable surface of genus $g$ having  
$n\geq n_g$ vertices are equivalent by means of 
diagonal transformations, and moreover this 
can be realized  among cubications without double edges. 
In particular, our main result is a consequence of the theorem of 
Nakamoto-Ota, in the case of closed orientable surfaces. 

\begin{corollary}
Let $C_0$ and $C_1$  be  marked cubications of a 
closed orientable surface 
$\Sigma$ that have the same invariants. Then there exists a sequence 
of flips relating them with the following properties: 
\begin{enumerate}
\item  start with a number of moves $b_1$ and $b_2$ such that 
$C_0$ and $C_1$ have the same number of vertices,  which 
is larger than $n_g$.  
\item use then diagonal transformations to connect 
the two cubications as in \cite{NO}. 
Decompose further the diagonal transformations into cubical flips as in the 
previous two pictures.  
\end{enumerate}
Thus, provided that the numbers of vertices of $C_i$ are large enough,  
there exists a sequence of flips that relating  them 
among cubications having no more than 
six extra vertices. In particular, the problem of finding  a connecting 
sequence of flips is algorithmically solvable.  
\end{corollary}

\begin{remark}
Although the methods of \cite{NO} do not work for 
non-orientable surfaces or surfaces with boundary, we expect 
that a result similar to the corollary above holds true in the general 
case. 
\end{remark}

  \begin{small} 
\bibliographystyle{plain}

\end{small}

 \end{document}